\documentclass[a4paper]{amsart}
 
\usepackage[utf8]{inputenc}
\usepackage[T1]{fontenc}

\usepackage{graphicx} 
\usepackage{caption}
\usepackage{subcaption}
\usepackage{booktabs} 
\usepackage{xcolor}
\usepackage{CJKutf8}
\usepackage{amsmath}
\usepackage{amsthm}
\usepackage{amssymb}
\usepackage{systeme}
\usepackage{hyperref}
\usepackage[top=3cm,bottom=3cm, left=2cm, right=2cm]{geometry}
\usepackage[squaren,Gray]{SIunits}
\usepackage{bbm}
\usepackage{stmaryrd}
\usepackage{pythontex}
\usepackage{listings}
\usepackage{setspace}
\usepackage{mathtools}
\usepackage{hhline}
\usepackage{mathrsfs}
\usepackage[overload]{empheq}

\newcommand{\R}{\mathbb{R}}

\newcommand{\dt}{\partial_t}

\newcommand{\dx}{\partial_x}
\newcommand{\dxx}{\partial_{xx}}
\newcommand{\dz}{\partial_z}
\newcommand{\dzz}{\partial_{zz}}

\newcommand{\s}{\sigma}

\title{A Paradigm for Well-Balanced Schemes for Traveling Waves Emerging in Parabolic Biological Models}

\author[M. Demircigil]{Mete Demircigil}
\address{Faculty of Science, Koç University, Istanbul, Turkey}
\email{mete.demircigil@math.cnrs.fr}

\author[B. Frabrèges]{Benoît Fabrèges}
\address{Institut Camille Jordan (ICJ), Universit\'e Claude Bernard Lyon~1, Lyon, France}
\email{fabreges@math.univ-lyon1.fr}

\date{}

\begin{document}
\maketitle
\begin{abstract}

We propose a methodology for designing well-balanced numerical schemes to investigate traveling waves in parabolic models from mathematical biology. 
We combine well-balanced techniques for parabolic models known in the literature with the so-called LeVeque-Yee formula as a dynamic estimate for the spreading speed. 
This latter formula is used to consider the evolution problem in a moving frame at each time step, where the equations admit stationary solutions, for which well-balanced techniques are suitable. 
Then, the solution is shifted back to the stationary frame in a well-balanced manner. 
We illustrate this methodology on parabolic reaction-diffusion equations, such as the Fisher/Kolmogorov-Petrovsky-Piskunov Equation, and a class of equations with a cubic reaction term that exhibit a transition from pulled to pushed waves.
We show that the numerical schemes capture in a consistent way simultaneously the wave speed and, to an extent, the so-called Bramson delay.
\end{abstract}

\section{Introduction}

Many biological systems may give rise to spatial spreading phenomena, such as for instance the invasion of rodents \cite{skellam1991} or the collective displacement of a population of \textit{Escherichia coli} cells in a micro-channel \cite{adler1966}. 
Mathematical modelling of this phenomena has contributed to a better understanding of these spreading phenomena and gives also rise to a quantitative framework, which may for instance account for the spreading speed. 
On the mathematical side, we may cite as models, that have attracted much attention, reaction-diffusion equations (\textit{e.g.} \cite{fisher1937,kolmogorov1937,aronson1975,hadeler1975,fife1977}) and parabolic models of chemotaxis (\textit{e.g.} \cite{keller1971,saragosti2011}). 
The biological question of the spreading of a population translates then into the mathematical question of the existence of a traveling wave, \textit{i.e.} solutions that are stationary in a moving frame of reference.

The mathematical investigation of traveling waves for biological models naturally raises the question of the numerical investigation of these traveling wave solutions. 
Numerically, this question is difficult, as it requires an accurate scheme over large domains in order to capture precisely the behavior of the solution over large time. 
Moreover, the numerical scheme should meet the following specifications: (i) find the accurate spreading speed; (ii) find the accurate traveling wave profile. 
In fact, these two questions are related, since mathematically the wave speed is intimately connected to the wave profile: wave profiles tend to decrease exponentially at their leading edge and their exponential decay parameter can be tied to the wave speed through a dispersion relation. 
Hence, we wish to find numerical schemes that capture at the same time accurately a wave speed, as well as an exponential decay, which requires great precision even at very low orders of magnitude, because of the nature of exponential decays.

In this article, we will use so-called well-balanced (WB) schemes in order to investigate these questions. 
The notion of WB schemes has been introduced by \cite{greenberg1996} and consists roughly speaking in numerical schemes that preserve the steady state of the continuous model, \textit{i.e.} the steady state of the numerical scheme matches the discretization of the steady state of the continuous model. Because of this property, these schemes compute faithfully the stationary profiles.
Many different approaches to such schemes exist and we refer the reader for instance to \cite{gosse2013} for a general discussion on these techniques. 
However in their essence, WB schemes accurately balance the flux terms and the source terms.
WB schemes have been investigated for parabolic equations \cite{gosse2017a,gosse2016}.
In these works, the strategy consists in interpolating the function in each cell with $\mathscr{L}$-splines (see \cite{gosse2018} for an overview on the concept of $\mathscr{L}$- splines), which correspond to the stationary solutions of the problem. 
The interpolation is then continous inside each cell, but at the nodes of the grid it may have $C^1$-discontinuities. 
The time integration then is computed by using the $C^1$ defect at each node. 
Of note, most referenced models are one-dimensional models, but in recent years the leap to higher dimensions has attracted interest \cite{despres2014,gosse2014,bretti2021}.

Yet, in the aformentioned studies, most often the schemes were designed to describe a stationary state, not a traveling wave. 
Whilst in some cases the propagation speed can be computed explicitly and thus be used in the implementation of the scheme, in most cases the propagation speed is not known \textit{a priori}. 
This adds to the difficulty of designing a WB scheme for traveling waves. 
For instance in the work \cite{calvez2017}, the authors have proposed a numerical scheme, which is well-balanced for a wave with speed $\s=0$, and have used it in order to capture traveling waves with nonzero speed.
Hence, the proposed scheme is \textit{stricto sensu} not WB for the observed traveling wave.

In the present article, in order to overcome this difficulty, we implement a strategy, which measures dynamically the spreading speed. We use an estimate of this spreading speed, which has been proposed in \cite{leveque1990} and has recently been refered to as LeVeque-Yee formula \cite{mascia2020}. This formula is based on the following observation: suppose that $u(t,x)=u(x-\s t)$ is a traveling wave profile, which admits at $\pm \infty$ fixed limit values $u(\pm \infty)$, then $u$ satisfies:
\begin{align*}
\dt u+\s\dx u=0.
\end{align*}
Integrating over $\R$ and rearranging the terms, we obtain the following identity:
\begin{align}
\label{LY-cont}
\s = -\frac{\int_\R \dt u dx}{u(+\infty)-u(\infty)}
\end{align}
The LeVeque-Yee formula is then merely a discretized version of Identity (\ref{LY-cont}):
\begin{align}
\label{LY}
\hat{\s}^n_{LY} = \frac{\Delta x}{t^n-t^{n-1}} \frac{\sum_i \left(u^{n-1}_i-u^n_i\right)}{u^n_I-u^n_0},
\end{align}
where $u$ is computed on a Cartesian grid $(t^n)_{n\in\mathbb{N}}\times (i\Delta x)_{i=0,\ldots,I}$. For solutions that are not traveling waves, the right hand-side of Identity (\ref{LY-cont}) remains well defined and can be interpreted as an estimate of their spreading speed.

The contribution of this article consists in combining the LeVeque-Yee formula as a dynamic speed estimate with the methodology of WB schemes for parabolic equations. The schemes presented in this article share the following common structure:
\begin{enumerate}
\item Estimate at each time step the spreading speed $\hat{\s}^n$ via the LeVeque-Yee Formula (\ref{LY}).
\item Consider over the time step the problem in the frame of reference moving at speed $\hat{\s}^n$ and use a WB scheme to integrate numerically the solution over the time step.
\item Given the solution integrated in time in the moving frame, shift back the solution to the stationary frame. To do so in a WB manner, interpolate in each cell the stationary solution, with the boundary conditions prescribed by the values obtained at the preceding step. Finally, in order to obtain the values on the grid of the stationary frame, use the value of the obtained stationary solution at the corresponding position. 
\end{enumerate}


In order to illustrate this methodology, we apply it to reaction-diffusion equations of the type:
\begin{equation}
    \label{eq.rd}
    \dt u-\dxx u=f(u).
\end{equation}
Many theoretical results are known about these types of equations, which will be used as a base of comparison for our numerical results.
First, we study the case $f(u)=u(1-u)$ of the Fisher/Kolmogorov-Petrovsky-Piskunov (F/KPP) equation \cite{fisher1937,kolmogorov1937}, which has been a prototype of traveling wave phenomena in parabolic equations. Then, we consider a cubic reaction function of the form $f(u)=u(1-u)(1+au)$ for a parameter $a>0$ introduced in \cite{hadeler1975}. This class is relevant, because it exhibits a transition from pulled to pushed waves (see Subsection \ref{ss.pulled.pushed}), which illustrates the behaviors of our numerical scheme on these regimes. Finally, we stress that this methodology is not restricted to reaction-diffusion equations and can be applied to parabolic equations with an advection term, such as for instance models of bacterial chemotaxis (\textit{e.g.} \cite{keller1971,saragosti2011}).

\subsection{Fisher/Kolmogorov-Petrovsky-Piskunov Equation}
\label{Sect-KPP}

The Fisher/Kolmogorov-Petrovsky-Piskunov (F/KPP) equation
\begin{align}
\label{fkpp-eqn}
\dt u - \dxx u = u(1-u),
\end{align}
introduced independently in \cite{fisher1937} and \cite{kolmogorov1937}, is a prototype of reaction-diffusion equations. It describes a population $u(t,x)$ of individuals undergoing unbiased motion, modeled through a diffusion operator, as well as competition among the population (\textit{e.g.} for ressources) through a logistic growth term $u(1-u)$: the higher the density $u$ becomes, the lower the growth rate $1-u$ becomes, with a saturation when the population reaches the maximal density $u=1$. It admits two steady states $u=0$ and $u=1$: the first one is unstable, whilst the second is stable. 

Interestingly, for velocity $\s$ Equation (\ref{fkpp-eqn}) admits nonnegative traveling wave solutions of the form:
\begin{align}
u_\s(-\infty)=1, u_\s(+\infty)=0\text{ and } u_\s(t,x)=u_\s(x-\s t)
\end{align}
if and only if $\s\geq \s^*=\s_{F/KPP}:=2$ (see \cite{fisher1937,kolmogorov1937,aronson1975}). These solutions are in fact invariant by translation and we will fix a specific solution by the convention that $u_\s(0)=\frac{1}{2}$. Hence the equation describes a linear spatial invasion of the state $u=0$ by the state $u=1$. Furthermore for a wide variety of initial data $u^0$, one can show that in a certain sense, which we will specify, $u$ the solution of Equation (\ref{fkpp-eqn}) converges to the traveling wave with minimal velocity $\s^*=\s_{F/KPP}$. Take for example an initial datum $u^0$, which satisfies:
\begin{align}
\label{s.int.DI}
\liminf_{x\to-\infty}u^0(x) > 0 \text{ and } \limsup_{x\to+\infty}
\frac{\dx u^0 (x)}{u^0(x)}\leq -1.
\end{align}
Then, define the position of the level set $x_c(t)$ for $c\in (0,1)$:
\begin{align}
x_c(t) = \sup \{ x | u(t,x)=c\}
\end{align}
Through a probabilistic interpretation of Equation (\ref{fkpp-eqn}) (see \cite{hamel2013,nolen2017} for a proof based on PDE argument), first Bramson \cite{bramson1978}, and later Uchiyama \cite{uchiyama1980} for inital conditions satisyfing (\ref{s.int.DI}), proved that there exists a constant $x_\infty$ depending on the initial datum and $c$, such that:
\begin{align}
\label{bramson}
x_c(t)=2t - \frac{3}{2}\ln(t)+x_\infty+o(1) \text{ and } \lim_{t\to +\infty}u(t,x+x_c(t)) = u_{\s^*}(x) .
\end{align}
The logarithmic delay in the asymptotic expansion of $x_c(t)$ has ever since been referred to as the Bramson delay. 
The result is in a sense very surprising, because it shows convergence to the traveling wave with speed $\s_{F/KPP}$ but in the frame shifted by a logarithmic correction term. 
As it has been observed in \cite{hamel2013}, this logarithmic correction term is due to subtle interactions between the diffusion and the reaction term, when $u\sim 0$. 
Hence, from a numerical point of view, capturing the Bramson delay is likely a very delicate task, since it requires an accurate computation of the profile of $u$ at very low orders of magnitude. 
As far as we know, we are not aware of any work on a numerical scheme for Equation (\ref{fkpp-eqn}), which addresses the question of the Bramson delay. 

\subsection{Reaction-Diffusion Equations with a Transition from Pulled to Pushed Waves}
\label{ss.pulled.pushed}

In \cite{hadeler1975}, the authors have introduced the reaction-diffusion equations with a cubic term:
\begin{align}
\label{semi-fkpp-eqn}
\dt u - \dxx u = u(1-u)(1+au).
\end{align}
For this model, the minimal wave speed solutions can be computed explicitly and it can be shown the following formula for the wave speed:
\begin{equation*}
\sigma^* = \left\{
\begin{array}{cc}
    2 & \text{ if }a\leq 2 \\
    \sqrt{\frac{a}{2}} + \sqrt{\frac{2}{a}} & \text{ if }a>2
\end{array}
\right. .
\end{equation*}
The dichotomy given by this formula extends into a dichotomy on the qualitative behavior of the waves, which is referred to as \textit{pulled} or \textit{pushed} regime \cite{stokes1976}. For $a\leq 2$, the wave is in the pulled regime: qualitatively the spreading is driven by the dynamics at the very leading edge of the front, where $u \sim 0$. For $a>2$, the wave is in the pushed regime: qualitatively, the spreading is subject to a contribution from the whole population.

These regimes are also reflected in the asymptotic expansion of the level sets $x_c(t)$. For $a>2$, as a consequence of the exponential convergence to the traveling wave proven in \cite{rothe1981}, we have:
\begin{align}
\label{bramson1}
x_c(t)=\sigma^* t +x_\infty+o\left(e^{-\omega t}\right),
\end{align}
for a constant $\omega >0$. For $a<2$, the level sets are again corrected by the logarithmic Bramson delay \cite{hamel2013}:
\begin{align}
\label{bramson2}
x_c(t)=2 t -\frac{3}{2}\ln(t)+x_\infty+o(1).
\end{align}
The critical case $a=2$ has recently been investigated in \cite{giletti2021,derrida2022} and following the terminology introduced in \cite{an2021} is referred to as the \emph{pushmi-pullyu} regime. The level sets are corrected by a logarighmic delay with a constant $-\frac{1}{2}$:
\begin{align}
\label{bramson3}
x_c(t)=2 t -\frac{1}{2}\ln(t)+x_\infty+o(1).
\end{align}

\subsection{Outline of the paper}

In this article, we combine the LeVeque-Yee formula with the methodology of WB schemes for 
parabolic equations in order to design numerical schemes that are WB for traveling waves. 
In Section \ref{sect-2}, we present our methodology for WB numerical schemes applied to reaction-diffusion equations.
In Section \ref{sect-3}, we numerically assess the WB schemes by applying them to the F/KPP Equation (\ref{fkpp-eqn}) and then to the cubic case (\ref{semi-fkpp-eqn}). 
The scheme is compared to an operator splitting approach. 
We show that our scheme captures accurately the propagation speed and we compare the asymptotic spreading with the logarithmic Bramson correction (\ref{bramson},\ref{bramson2},\ref{bramson3}) to the asymptotic expansion in the pulled case. As a result we will show that the presented methodology works particularily well in the pulled case, but does not outperform standard schemes in the pushed case.

\section{A WB Scheme for Reaction-Diffusion Equations}
\label{sect-2}

In order to solve numerically reaction-diffusion equations (\ref{eq.rd}) in a WB manner, we propose a numerical scheme based on \cite{gosse2017a} with the difference that we also investigate an implicit Euler method for time integration instead of the standard explicit Euler method.
In fact, as it has been pointed out in \cite{gosse2017a}, for the explicit Euler method, one requires essentially a parabolic CFL condition in order to guarantee stability. 
For the implicit Euler method, through numerical investigation, we observe that a hyperbolic CFL condition is sufficient. 
Moreover, resolving the time integration implicitly comes computationally at not too high of a cost, because of the special tridiagonal structure of the problem.

We consider a Cartesian grid $(t^n)\times(x_i)$, with $x_i=i\Delta x$, for $i \in \mathbb{Z}$ for a given $\Delta x$. 
The time points $(t^n)$ are set dynamically because of the lack of an \textit{a priori} bound to satisfy a CFL condition. We set $\Delta t^n=t^{n+1}-t^n$. 
We introduce the parallel grid with points $z_i=i \Delta z$, $\Delta z=\Delta x$, $\bar{C}_{i+\frac{1}{2}}=\left(z_{i},z_{i+1}\right)$ and consider $u^n_i\approx u(t^n,x_i)$.

The procedure of the numerical scheme can be summarized as follows:
\begin{enumerate}
\item We use LeVeque-Yee Formula (\ref{LY}) on the profile $u$ in order to obtain an estimate $\hat{\s}^n$ of the propagation speed.
\item Given $\hat{\s}^n$, on the time interval $(t^n,t^{n+1})$ we consider $u$ in the moving frame $(t,z)=(t,x-\hat{\s}^n(t-t^n))$ and denote it by $\bar{u}$, which leads to Equation:
\begin{align}
\label{mv-fkpp}
\dt \bar u- \hat{\s}^n\dz \bar  u-\dzz \bar u = f(\bar u).
\end{align}
In order to deal with the nonlinear term in the right hand-side of Equation (\ref{eq.rd}), we freeze the non-linear contribution to the term on each mesh $(t^n,t^{n+1})\times \bar{C}_{i+\frac{1}{2}}$ by considering that $f(\bar u)\approx \bar u \bar{F^n}_{i+\frac{1}{2}}$, where $\bar{u}^n_{i+\frac{1}{2}} = \frac{{u}^n_i+{u}^n_{i+1}}{2}$, $\bar{F^n}_{i+\frac{1}{2}}=F(\bar{u}^n_{i+\frac{1}{2}})$ and  $F(u)=\frac{f(u)}{u}$. 
As noted in the Introduction, the dynamics of the F/KPP Equation (\ref{fkpp-eqn}) or Equation (\ref{semi-fkpp-eqn}) with $a\leq 2$ is predominantly set by the dynamics at the leading edge of the profile, where $u\sim 0$. 
There, the linearization of Equation (\ref{eq.rd}) approximates extremely well the nonlinear Equation. 
Hence it is reasonnable to expect that this way of handling numerically the nonlinearity will not influence much the dynamics of the propagation. 
Thanks to the freezing of the nonlinearity, we obtain a linear problem for which we can apply the methodology of $\mathscr{L}$-spline interpolation, as proposed in \cite{gosse2017a}, leading to the values $\bar{u}^{n+1}_i$, which approximate $\bar{u}(t^{n+1},z_i)$.
\item It remains to shift back to the stationary frame $(t,x)$. The values $\bar{u}^{n+1}_{i}$ are used to extrapolate a value for $u(t^{n+1},x_i)$, \textit{i.e.} $\bar{u}(t^{n+1},z_i+\hat{\s}^n\Delta t^n)$. To do so we consider the solution $\hat{u}$ of the following stationary problem in each cell $\bar{C}_{i+\frac{1}{2}}$:
\begin{align}[left = \empheqlbrace\,]
\label{stat-fkpp}
&-\hat{\s}^n\dz \hat{u}-\dzz \hat{u} = \hat{u}\hat{F^{n}}_{i+\frac{1}{2}}, \text{ for }z\in \bar{C}_{i+\frac{1}{2}} \\
&\hat{u}(z_i)=\bar{u}^{n+1}_i \text{ and }\hat{u}(z_{i+1})=\bar{u}^{n+1}_{i+1},\nonumber
\end{align}
with $\hat{u}^{n+1}_{i+\frac{1}{2}} = \frac{\bar{u}^{n+1}_i+\bar{u}^{n+1}_{i+1}}{2}$ and $\hat{F}^{n+1}_{i+\frac{1}{2}}=F\left(\hat{u}^{n+1}_{i+\frac{1}{2}} \right)$. Finally, we set:
\begin{align*}
u^{n+1}_i=\hat{u}(z_i+\hat{\s}^n\Delta t^n).
\end{align*}
\end{enumerate}
In order to avoid taking values from cells which would not be direct neighbors, we enforce the following CFL condition:
\begin{align}
\label{fkpp-cfl}
\hat{\s}^n\Delta t^n \leq \Delta x.
\end{align}

\subsection{LeVeque-Yee Formula for F/KPP Equation}
\label{LY-FKPP}

We apply the LeVeque-Yee formula to the profile $u$ in order to obtain an estimate for the propagation speed:
\begin{align}
\label{LY-fkpp}
\hat{\s}^n = \frac{\Delta x}{\Delta t^{n-1}} \frac{\sum_i \left(u^{n-1}_i-u^n_i\right)}{u^n_I-u^n_0}.
\end{align}
The leftmost (resp. rightmost) value $u^{n}_0$ (resp. $u^{n}_I$) will in practice (possibly after a short transitory phase) be very close to the steady state $1$ (resp. the steady state $0$). 
Moreover in order to satisfy CFL condition (\ref{fkpp-cfl}), we set dynamically the time points of the grid, in order to avoid imposing a restriction on the speed estimate.

\subsection{A WB Scheme in the Moving Frame}

As explained above, we consider the evolution of $f$ in the moving frame $(t,z)=(t,x-\hat{\s}^n\Delta t^n)$, \textit{i.e.} $\bar{u}(t,z)=u(t,z+\hat{\s}^n(t-t^n))$, where $\bar{u}$ solves Equation (\ref{mv-fkpp}). In order to approximate numerically this solution, we use the method of $\mathscr{L}$-spline interpolation proposed in \cite{gosse2017a}.
The scheme is based either on the explicit integration in time rule:
\begin{align}
\label{fkpp-update-explicit}
\bar{u}^{n+1}_i = \bar{u}^{n}-\frac{\Delta t^n}{\Delta x}\left(L^{n}_{i+\frac{1}{2}} - R^{n}_{i-\frac{1}{2}}\right),
\end{align}
or on the implicit integration in time rule:
\begin{align}
\label{fkpp-update}
\bar{u}^{n+1}_i = \bar{u}^{n}-\frac{\Delta t^n}{\Delta x}\left(
L^{n+1}_{i+\frac{1}{2}} - R^{n+1}_{i-\frac{1}{2}}
 \right),
\end{align}
where $R^{n+1}_{i-\frac{1}{2}}, L^{n+1}_{i+\frac{1}{2}}, R^{n}_{i-\frac{1}{2}}, L^{n}_{i+\frac{1}{2}}$ correspond to numerical fluxes at the point $z=z_i$ of two interpolated functions in $\bar{C}_{i-\frac{1}{2}}$ and $\bar{C}_{i+\frac{1}{2}}$. With the implicit time integration, the values $R^{n+1}_{i-\frac{1}{2}}, L^{n+1}_{i+\frac{1}{2}}$ are obtained by solving the linear stationary problem:
\begin{align}[left = \empheqlbrace\,]
\label{mv-stat-fkpp}
&-\hat{\s}^n\dz w-\dzz w = w\bar F^{n}_{i+\frac{1}{2}}, \text{ for }z\in \bar{C}_{i+\frac{1}{2}} \\
&w(z_i)=\bar{u}^{n+1}_i \text{ and }w(z_{i+1})=\bar{u}^{n+1}_{i+1}.
\nonumber 
\end{align}
Given the solution $w$, we then set $R^{n+1}_{i+\frac{1}{2}} := w'(z_{i+1}^-)$ and $L^{n+1}_{i+\frac{1}{2}}:= w'(z_i^+)$.
One can view $L^{n+1}_{i+\frac{1}{2}} - R^{n+1}_{i-\frac{1}{2}}$ as a measure of the defect of $C^1$-smoothness at $z_i$. 
The time integration (\ref{fkpp-update}) will therefore be stationary if $\bar{u}$ solves Equation (\ref{mv-stat-fkpp}) on the whole space, which in particular implies that it is $C^1$ as a solution to a second-order elliptic equation. Thus this procedure is well-balanced in the sense that the numerical scheme admits as stationary state the discretization of the continuous stationary state. For the explicit case, one merely changes the boundary conditions in (\ref{mv-stat-fkpp}) to $w(z_i)=\bar{u}^{n}_i$ and $w(z_{i+1})=\bar{u}^{n}_{i+1}$ and follows the same procedure.

Numerically, we have observed that the explicit time integration is stable under the parabolic CFL condition that:
\begin{align*}
    \Delta t^n \leq \frac{1}{2}\Delta x^2
\end{align*}
The coefficient $\frac{1}{2}$ can be slightly higher, but for $0.6$, we have observed instability. This empirical CFL condition is in the same spirit as an explicit CFL condition in \cite{gosse2017a}. For the implicit time integration, we have observed that under the hyperbolic CFL condition (\ref{fkpp-cfl}), the scheme remains stable. 

Let us now describe the implicit (linear) relation between $(L^{n+1}_{i+\frac{1}{2}},R^{n+1}_{i+\frac{1}{2}})$ and $(\bar{u}^{n+1}_i,\bar{u}^{n+1}_{i+1} )$. 
The same approach applies \textit{mutatis mutandis} for the explicit case. 
We solve Problem (\ref{mv-stat-fkpp}) with respect to $(\bar{u}^{n+1}_i,\bar{u}^{n+1}_{i+1} )$. 
On the cell $\bar{C}_{i+\frac{1}{2}}$, it is a simple second-order differential equation with constant coefficients. We compute the roots of its characteristic polynomial:
\begin{align}
\mu^2+\hat{\s}^n\mu+\bar F^n_{i+\frac{1}{2}}=0
\end{align}
The discriminant is $\Delta_{i+\frac{1}{2}}^n =\left(\hat{\s}^{n}\right)^2-4 \bar F^n_{i+\frac{1}{2}}$ and three cases exist:
\begin{enumerate}
\item Suppose $\Delta_{i+\frac{1}{2}}^n>0$, then set $ \mu^{n }_{i+\frac{1}{2},\pm} = \frac{-\hat{\s}^n\pm \sqrt{\Delta_{i+\frac{1}{2}}^n}}{2}$. Then for $z\in \bar{C}_{i+\frac{1}{2}}$, we have:
\begin{align*}
w(z) = a^- e^{\mu^{n }_{i+\frac{1}{2},-} z}+a^+e^{\mu^{n }_{i+\frac{1}{2},+}z}
\end{align*}
Therefore by setting the following matrices:
\begin{align}
\label{fkpp-mat}
M = 
\begin{pmatrix}
e^{\mu^{n }_{i+\frac{1}{2},-} z_i} & e^{\mu^{n }_{i+\frac{1}{2},+} z_i}\\
e^{\mu^{n }_{i+\frac{1}{2},-} z_{i+1}} & e^{\mu^{n }_{i+\frac{1}{2},+} z_{i+1}}
\end{pmatrix}
\text{ and }
\tilde{M} = 
\begin{pmatrix}
\mu^{n }_{i+\frac{1}{2},-}e^{\mu^{n -}_{i+\frac{1}{2}} z_i} & \mu^{n }_{i+\frac{1}{2},+}e^{\mu^{n +}_{i+\frac{1}{2}} z_i}\\
\mu^{n }_{i+\frac{1}{2},-}e^{\mu^{n -}_{i+\frac{1}{2}} z_{i+1}} & \mu^{n }_{i+\frac{1}{2},+}e^{\mu^{n +}_{i+\frac{1}{2}} z_{i+1}}
\end{pmatrix},
\end{align}
we obtain the following two systems:
\begin{align*}
\begin{pmatrix}
\bar{u}_i ^{n+1}\\
\bar{u}_{i+1}^{n+1}
\end{pmatrix}=
M
\begin{pmatrix}
a^-\\
a^+
\end{pmatrix}
\text{ and }
\begin{pmatrix}
L_{i+\frac{1}{2}} ^{n+1}\\
R_{i+\frac{1}{2}}^{n+1}
\end{pmatrix}=
\tilde{M}
\begin{pmatrix}
a^-\\
a^+
\end{pmatrix}.
\end{align*}
Hence, by setting $S_{i+\frac{1}{2}} = \tilde{M}M^{-1}$, we obtain the following identity:
\begin{align}
\label{implicit}
\begin{pmatrix}
L_{i+\frac{1}{2}} ^{n+1}\\
R_{i+\frac{1}{2}}^{n+1}
\end{pmatrix}
=
S_{i+\frac{1}{2}} \begin{pmatrix}
\bar{u}_i ^{n+1}\\
\bar{u}_{i+1}^{n+1}
\end{pmatrix}
\end{align}
\item Suppose $\Delta_{i+\frac{1}{2}}^n=0$, then we set $\mu_{i+\frac{1}{2}}^n:=-\frac{\hat{\s}^n}{2}$ and the solution is of the shape:
\begin{align*}
w(z) = (az+b) e^{\mu_{i+\frac{1}{2}}^n z}.
\end{align*}
The definition of the matrix $S_{i+\frac{1}{2}}$ follows \textit{mutatis mutandis} as precedingly.
\item In the case $\Delta_{i+\frac{1}{2}}^n<0$, the solution $w$ is of the shape:
\begin{align*}
w(z) = \left(a\cos\left( \sqrt{-\Delta_{i+\frac{1}{2}}^n}  z\right) +b\sin\left( \sqrt{-\Delta_{i+\frac{1}{2}}^n} z\right)\right) e^{-\frac{\hat{\s}^n z}{2}},
\end{align*}
and again the definition of the matrix $S_{i+\frac{1}{2}}$ then follows \textit{mutatis mutandis} as precedingly.
\end{enumerate}

Relation (\ref{implicit}) can be expressed explicitly, but we will not do so for the sake of concision. Then, combining Relation (\ref{implicit}) with the time integration (\ref{fkpp-update}) leads to a tridiagonal system, which can be inverted via Thomas' algorithm.

\subsection{A WB Shift to the Stationary Frame}

In the next step, we need to compute the values $u^{n+1}_i$ (in the stationary frame) from the values $\bar{u}^{n+1}_i$ (in the moving frame $(t,z)=(t,x-\hat{\s}^n(t-t ^n))$). To do so, we solve Equation (\ref{stat-fkpp}), which we recall:
\begin{align*}[left = \empheqlbrace\,]
&-\hat{\s}^n\dz \hat{u}-\dzz \hat{u} = \hat{u}\hat F^{n+1}_{i+\frac{1}{2}}, \text{ for }z\in \bar{C}_{i+\frac{1}{2}} \\
&\hat{u}(z_i)=\bar{u}^{n+1}_i \text{ and }\hat{u}(z_{i+1})=\bar{u}^{n+1}_{i+1},
\end{align*}
where we recall that $\hat F^{n+1}_{i+\frac{1}{2}} = F\left( \hat u^{n+1}_{i+\frac{1}{2}}\right)$ is an updated version of the nonlinearity with $\hat{u}^{n+1}_{i+\frac{1}{2}}=\frac{\bar u^{n+1}_{i}+ \bar u  ^{n+1}_{i+1}}{2}$. Then, we set:
\begin{align}
u^{n+1}_i=\hat{u}(z_i+\hat{\s}^n\Delta t^n).
\end{align} 
Enforcing the CFL condition (\ref{fkpp-cfl}) avoids taking values from cells which would not be direct neighbors. Equation (\ref{stat-fkpp}) can be solved in the same manner as previously. One takes $N_{i+\frac{1}{2}}=M_{i+\frac{1}{2}}$ and for $(\phi_1,\phi_2)$ the fundamental system of solutions, that we have used above:
\begin{align}
\tilde{N}_{i+\frac{1}{2}} =
\begin{pmatrix}
\phi_1(z_i+\hat{\s}^n\Delta t^n)  & \phi_2(z_i+\hat{\s}^n\Delta t^n) 
\end{pmatrix}
\end{align}
This then leads to the following $1\times 2$ matrix:
\begin{align*}
T_{i+\frac{1}{2}}:=  \tilde{N}_{i+\frac{1}{2}}N_{i+\frac{1}{2}}^{-1},
\end{align*}
and:
\begin{align}
u^{n+1}_i=
T_{i+\frac{1}{2}}
\begin{pmatrix}
\bar{u}^{n+1}_i \\
\bar{u}^{n+1}_{i+1}
\end{pmatrix}.
\end{align}

\section{Numerical Assessements}
\label{sect-3}

We start by comparing the preceding WB scheme for F/KPP Equation (\ref{fkpp-eqn}) with two other numerical schemes, which are described just below. A similar comparison is carried out for the Equation (\ref{semi-fkpp-eqn}).

The code for the numerical schemes in this Section can be found in \cite{fabregesSchemePLMLAB}.

\subsection{Alternative Schemes}

\subsubsection{OS Scheme with Crank-Nicolson Method}
\label{sect-fkpp-os}

The first scheme is based on an operator splitting (Strang splitting) of the heat operator and the reaction term:
\begin{align*}
    u^{n+1}=R_{\frac{\Delta t}{2}}e^{-\Delta t \dxx }R_{\frac{\Delta t}{2}}u^n
    \end{align*}
For the reaction term, in general we may use an implicit Euler integration, but in the special case for the F/KPP equation, we can get an exact inegration as follows:
\begin{align*}
    R_{\frac{\Delta t}{2}}v=\frac{e^{\frac{\Delta t}{2}}}{e^{\frac{\Delta t}{2}} - 1 + \frac{1}{v}}.
\end{align*}
The heat operator is integrated via the Crank-Nicolson method with $v^{n+1}=e^{-\Delta t \dxx }v^n$ satisfying for all $i$:
\begin{align*}
\frac{v^{n+1}_i-v^n_i}{\Delta t} = \frac{v^{n+1}_{i-1}-2v^{n+1}_i+v^{n+1}_{i+1}}{2\Delta x^2}+ \frac{v^{n}_{i-1}-2v^{n}_i+v^{n}_{i+1}}{2\Delta x^2}.
\end{align*}

\subsubsection{WB Scheme for the "0-Wave"}
\label{sect-fkpp-wb-0wave}

The second scheme is based on the same method as the WB scheme mentioned before, but we impose that $\hat{\s}^n=0$. This choice corresponds to the choice that has been investigated in \cite{gosse2017a}.
As mentioned in the Introduction, this approach is only WB for a stationary state, or equivalently a traveling wave with speed $0$. For the implementation of this scheme, one does not need to shift from the moving to the stationary frame, as both coincide. Or, with the notations from above, we have:
\begin{align*}
T_{i+\frac{1}{2}}=  \begin{pmatrix}
1 \\
0
\end{pmatrix}.
\end{align*}

\subsection{The Asymptotic Propagation Speed}

The F/KPP Equation (\ref{fkpp-eqn}) has an asymptotic progation speed $\sigma$ of $2$. We compare the ability of different schemes to capture this asymptotic velocity: the WB scheme of Section \ref{sect-2} in its explicit and implicit version, the OS scheme of Section \ref{sect-fkpp-os} and the WB scheme for the "0-Wave" of Section \ref{sect-fkpp-wb-0wave} also in its explicit and implicit version.

In order to capture the asymptotic propagation speed, the F/KPP equation is solved over the time interval $[0, 1500]$. Because the expected spreading speed is $\s=2$, the domain considered should be at least twice bigger: here, we take the domain $[0, 3080]$. 
This large domain avoids roughly speaking boundary effects. The initial condition is the sigmoid function $u_0(x) = 1 - \frac{1}{1 + e^{-3(x - 40)}}$.

\begin{figure}
\begin{center}
\includegraphics[width=12cm]{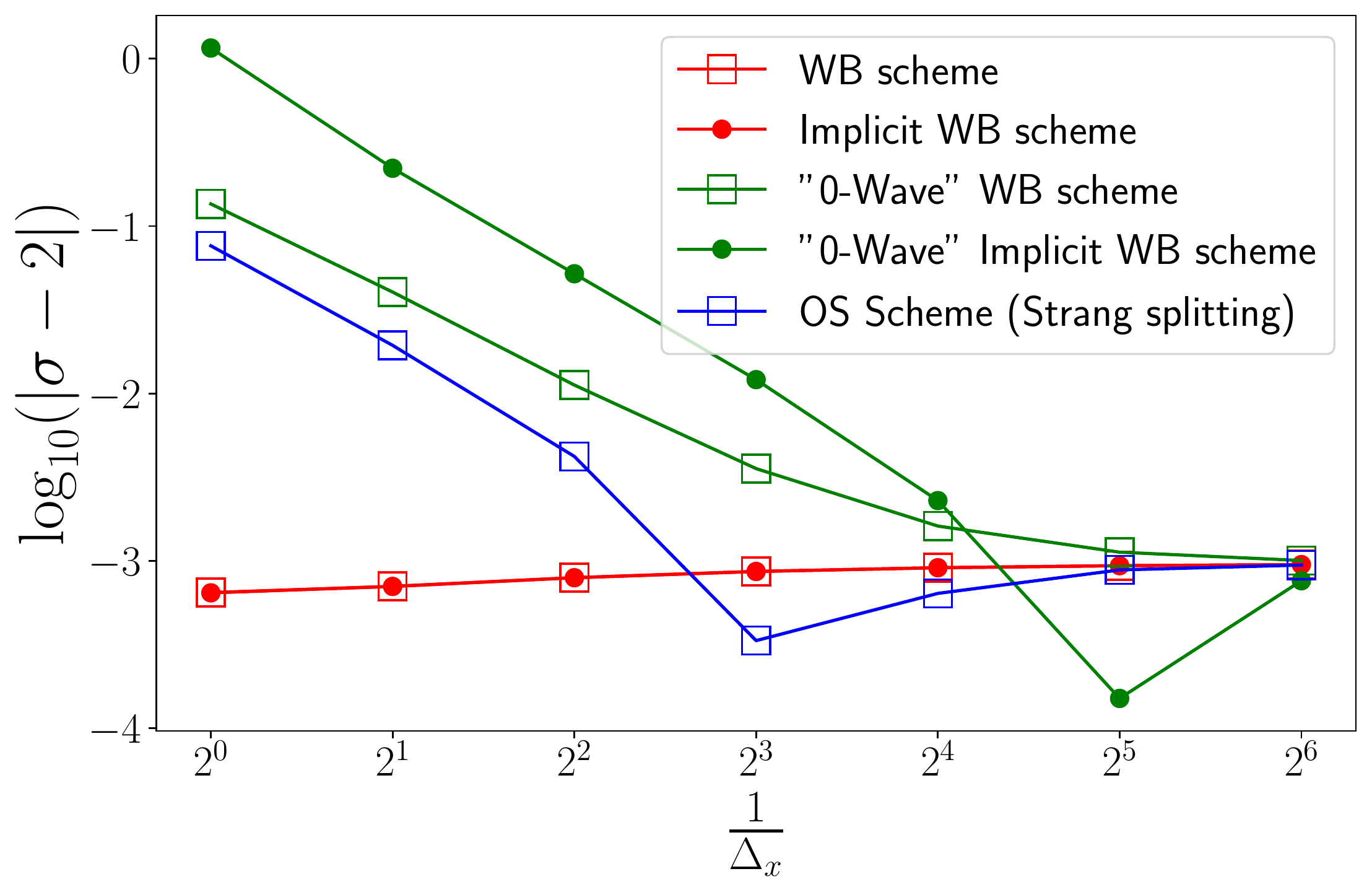}
\caption{Evolution of LeVeque-Yee velocity $\sigma$ for the F/KPP equation at time $t=1500$ for different schemes with respect to the space discretization size $\Delta_x$. \label{Fig-FKPP-Asymptotic-Speed} }
\end{center}
\end{figure}

Figure \ref{Fig-FKPP-Asymptotic-Speed} represents the velocity $\hat{\s}^n$ obtained with the LeVeque-Yee formula at the final time $t=1500$ for the schemes considered here.
With a coarse mesh, the asymptotic propagation speed obtained is far from $2$ for some schemes and the spatial domain had to be extended to $[0, 6080]$ in order to avoid boundary effects for these schemes.
The same time step is chosen for all the considered schemes at each iteration. It is the one satisfying all the CFL conditions which is the parabolic one, $\Delta_t^n \leq \frac{1}{2}\Delta_x^2 $, needed for the explicit WB schemes.
As the meshes become finer, all the numerical schemes seem to converge to the same result, which is extremely close to the asymptotic propagation speed of $2$. Yet, one can observe that the speed does not quite converge to $2$. Heuristically, this discrepancy may be a consequence of the Bramson delay. Indeed, differentiating formally the asymptotic expansion (\ref{bramson}) leads $\dot{x_c}(t)=2-\frac{3}{2t}+o\left(\frac{1}{t}\right)$. Then plugging in the value $t=1500$ would explain the observed discrepancy of the order $10^{-3}$. 
Moreover, it can be noted that the WB scheme, both in its explicit and implicit version, captures extremely well even for very coarse meshes the asymptotic behavior that is captured by the other schemes for much finer meshes.

\subsection{Capturing the Bramson Delay}

We recall that the Bramson delay is the logarithmic correction that appears in the expression of the position of a level set $x_c(t)$ for $c \in (0, 1)$ given by Equation (\ref{bramson}).
In the simulations, we track the position of the level set $x_{0.5}(t)$ by linearly interpolating the solution in the mesh cell where the level set is. The coefficient $\frac{3}{2}$ of the Bramson delay is approximated by fitting the model $m_{\alpha, \beta, \gamma}(t) = \alpha \log(t) + \beta + \frac{\gamma}{\sqrt{t}}$ to $x_{0.5}(t) - 2t$ over the last half of the simulation. As the next term in the asymptotic expansion (\ref{bramson}) is $-\frac{3\sqrt{\pi}}{\sqrt{t}}$ \cite{nolen2019}, the $\frac{\gamma}{\sqrt{t}}$ is added in order to improve the fitting of the model.

\begin{figure}
\begin{center}
\includegraphics[width=12cm]{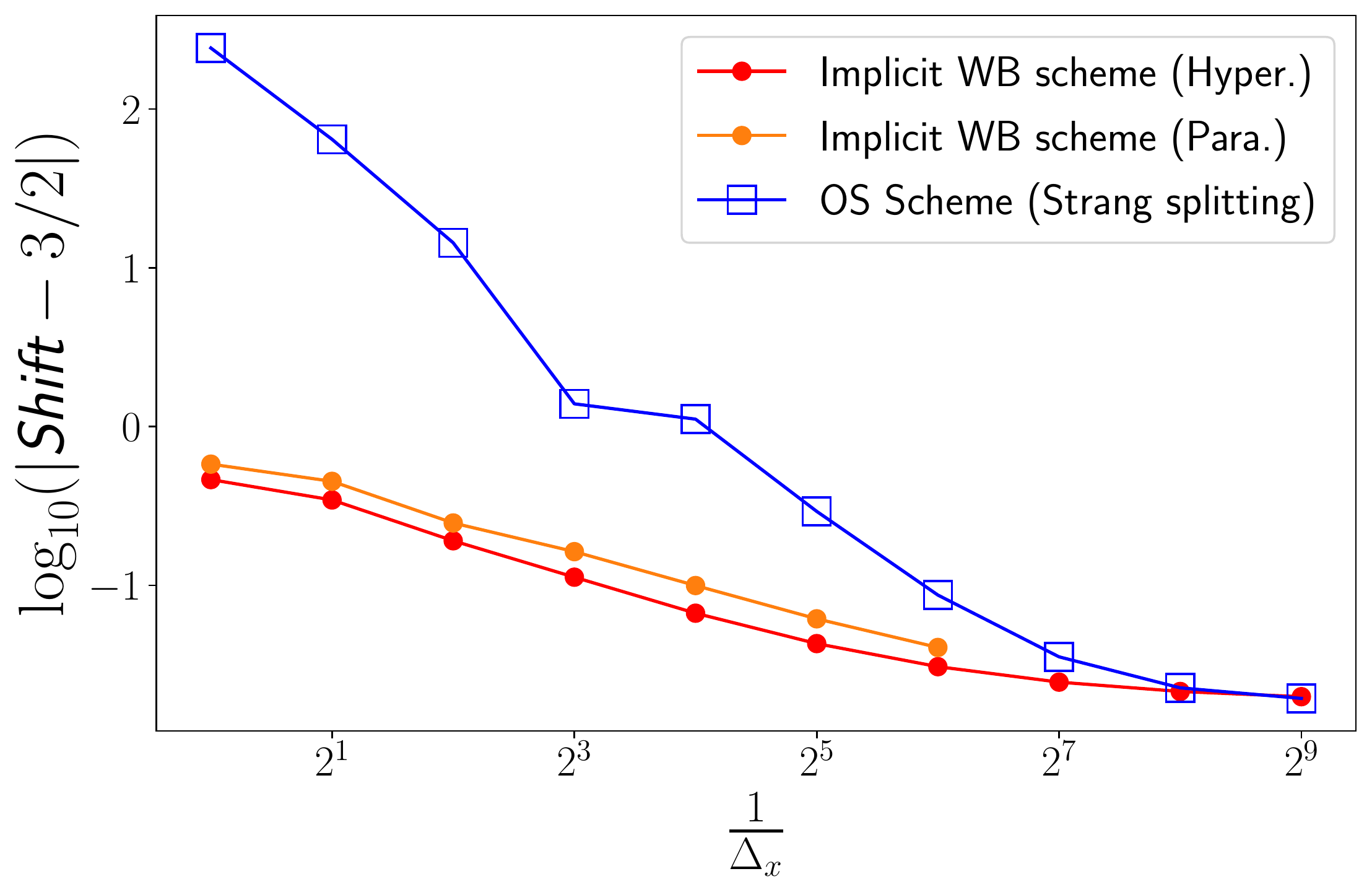}
\caption{Fitting of the coefficient $3/2$ of the Bramson delay with the implicit WB scheme and the OS scheme for the F/KPP Equation. \label{Fig-FKPP-Bramson-Shift} }
\end{center}
\end{figure}

We try to capture this Bramson delay by refining the discretization of the spatial domain.
Using schemes that depend on a parabolic CFL condition results in a huge number of iterations in time which becomes rapidly intractable in a reasonable amount of time.
Moreover, the "0-Wave" implicit WB scheme does not capture the asymptotic propagation speed as well as the implicit WB scheme.
Thus, we only fit the Bramson delay for the implicit WB scheme and the OS scheme (see Figure \ref{Fig-FKPP-Bramson-Shift}).
The WB scheme and the OS scheme are of order $1$ and $2$ in time respectively. 
In order to compare these schemes, we also plot the fitting of the WB scheme using small enough time step to mimic a second order in time scheme.
In this case, the time step is $\Delta_t^n = \min(0.05, \frac{1}{2} \Delta_x^2)$. Yet, this parabolic CFL condition again leads to computational costs that are too high, when $\Delta_x$ is too small. Thus we have computed the implicit WB scheme with parabolic CFL condition only for $\Delta_x\geq 2^{-6}$. 
For all computed schemes, the shift obtained seems to converge to a value slightly different from $3/2$. One possible explanation for this discrepancy might be higher order terms in the expansion (\ref{bramson}) (see for instance \cite{graham2019}).
Here again, the implicit WB scheme captures the Bramson delay even for coarse meshes.
Moreover, for the implicit WB scheme, there is no significant difference between a time step $\Delta_t^n = \mathcal{O}(\Delta_x)$ and a time step $\Delta_t^n = \mathcal{O}(\Delta_x^2)$.

\subsection{Asymptotic Propagation Speed and Bramson Delay for the Cubic Case}

For Equation (\ref{semi-fkpp-eqn}), we compare the solutions obtained with the WB and the OS schemes in these three different regimes
\begin{itemize}
    \item pulled regime with $a = 1$,
    \item pushmi-pullyu transition regime with $a = 2$,
    \item pushed regime with $a = 3$.
\end{itemize}


\begin{figure}
\begin{center}
\includegraphics[width=12cm]{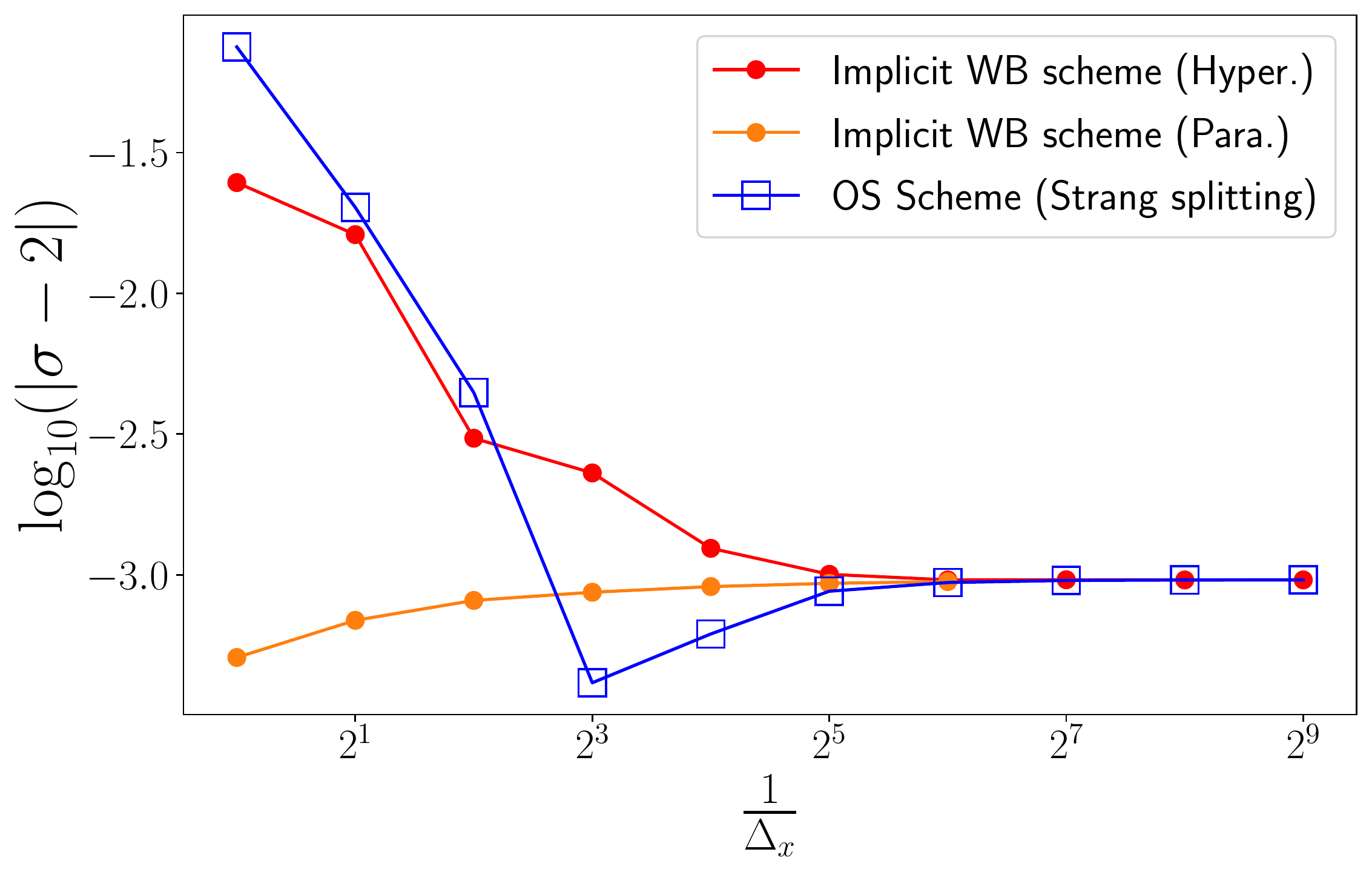}
\caption{Asymptotic propagation speed for the cubic case in a pulled regime ($a = 1$). For the OS scheme, we observe at $\Delta x=2^{-3}$ a change in sign of $\s-2$. \label{Fig-FKPP-Asymptotic-Speed-1+u} }
\end{center}
\end{figure}

\begin{figure}
\begin{center}
\includegraphics[width=12cm]{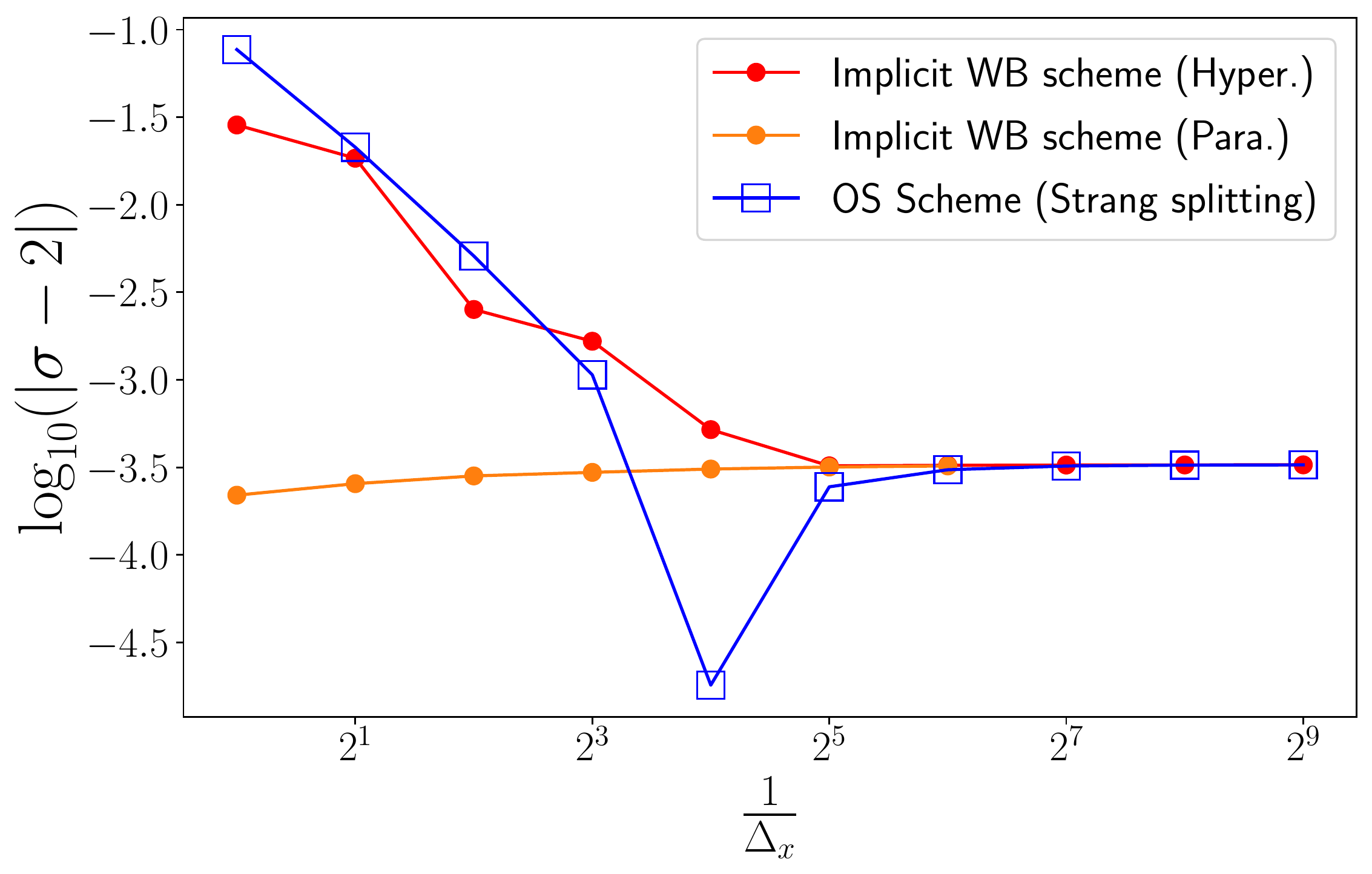}
\caption{Asymptotic propagation speed for the cubic case in the pushmi-pullyu transition regime ($a = 2$). For the OS scheme, we observe at $\Delta x=2^{-4}$ a change in sign of $\s-2$. \label{Fig-FKPP-Asymptotic-Speed-1+2u} }
\end{center}
\end{figure}

\begin{figure}
\begin{center}
\includegraphics[width=12cm]{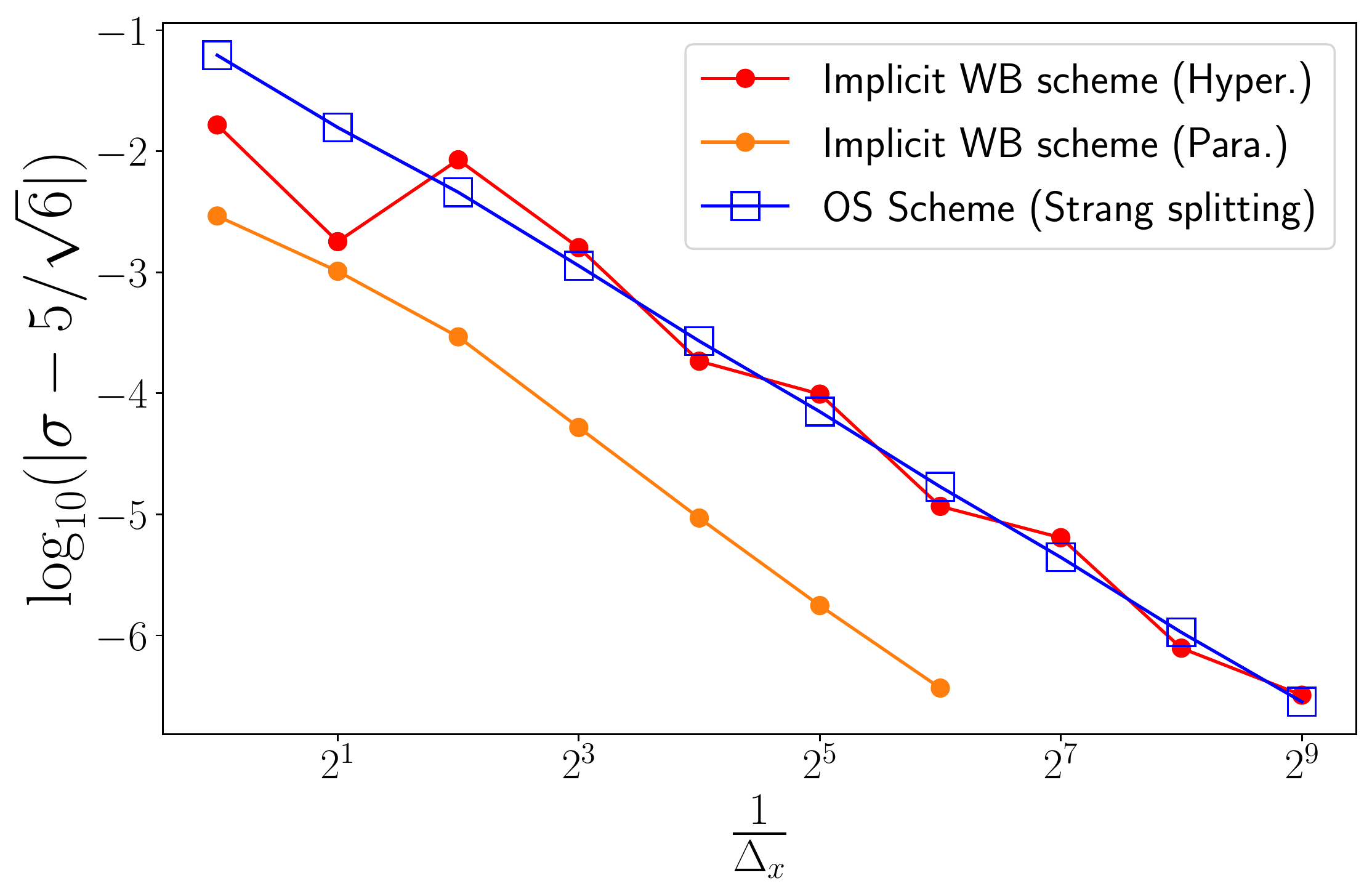}
\caption{Asymptotic propagation speed for the cubic case in a pushed regime ($a = 3$). \label{Fig-FKPP-Asymptotic-Speed-1+3u} }
\end{center}
\end{figure}

Figures \ref{Fig-FKPP-Asymptotic-Speed-1+u}, \ref{Fig-FKPP-Asymptotic-Speed-1+2u}, and \ref{Fig-FKPP-Asymptotic-Speed-1+3u}, illustrate the convergence of the numerical asymptotic speed in the three regimes.
In all cases, the velocities are well captured by the schemes, but the implicit WB scheme with a parabolic CFL performs especially better than the other ones. However, as already mentioned, the computational cost of the parabolic CFL is too high to use it with fine meshes.
In the pushed regime ($a=3$), it can be noted that for all schemes as the meshes get finer, the speed gets closer and closer to the predicted speed (Figure \ref{Fig-FKPP-Asymptotic-Speed-1+3u}). This is likely due to the exponential convergence in time given by asymptotic expansion (\ref{bramson1}), which states that the expansion speed at $t=1500$ of the continuous solution should be extremely close to the predicted speed, contrary to what happens in the pulled case due to the Bramson delay. Moreover, both schemes under hyperbolic CFL condition perform roughly equally well in the pushed case.

\begin{figure}
\begin{center}
\includegraphics[width=12cm]{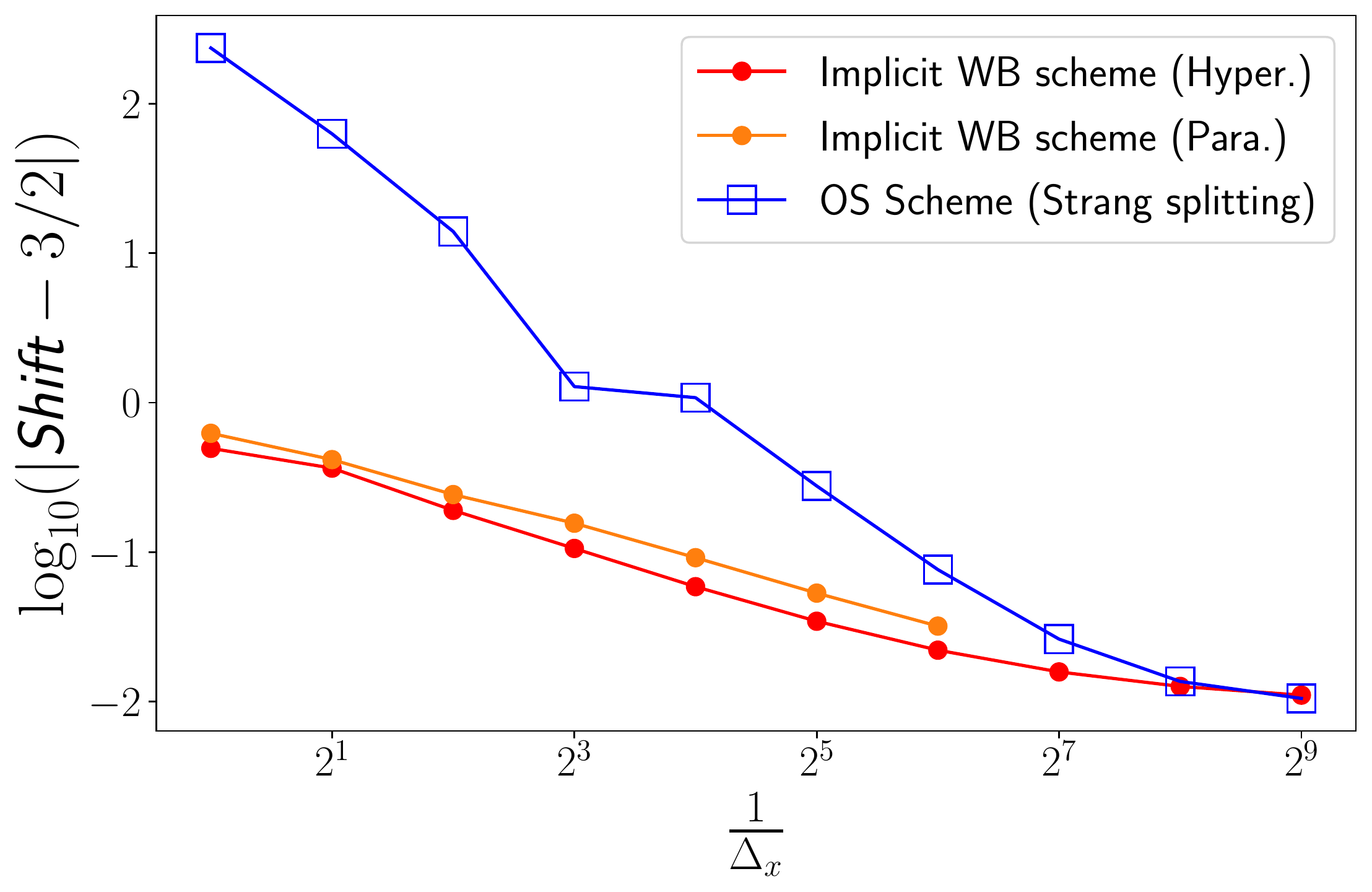}
\caption{Bramson delay for the cubic case in a pulled regime ($a = 1$). \label{Fig-FKPP-Bramson-Shift-1+u} }
\end{center}
\end{figure}

\begin{figure}
\begin{center}
\includegraphics[width=12cm]{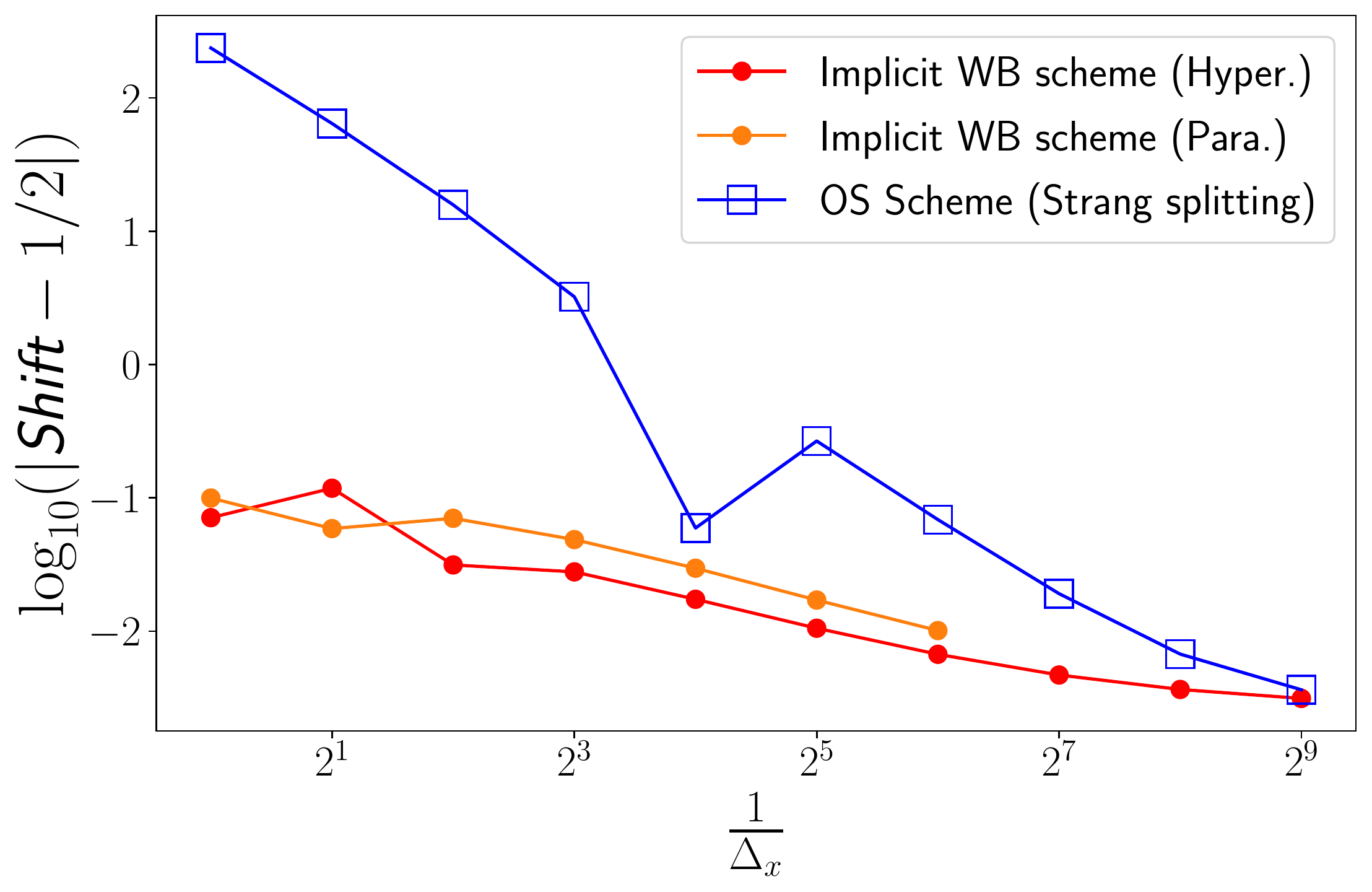}
\caption{Bramson delay for the cubic case in the pushmi-pullyu transition regime ($a = 2$). \label{Fig-FKPP-Bramson-Shift-1+2u} }
\end{center}
\end{figure}

Concerning the Bramson delay, the results for the pulled regime (see Figure \ref{Fig-FKPP-Bramson-Shift-1+u}) are very similar to the ones for the F/KPP equation. In particular, the WB scheme approximates well the shift even for coarse meshes. The Bramson delay of $1/2$ in the pushmi-pullyu transition regime (see Figure \ref{Fig-FKPP-Bramson-Shift-1+2u}) is captured by all schemes. Yet again, the WB scheme recover the $1/2$ shift even for coarse meshes, whilst the OS scheme catches up only for finer meshes.

\section*{Acknowledgements}

The authors are extremely grateful to Vincent Calvez and Thierry Dumont, who have been essential to the present work.\\
{This project has received funding from the European Research Council (ERC) under the European Union’s Horizon 2020 research and innovation programme (grant agreement No 865711).}


\bibliographystyle{siam}
\bibliography{biblio}

\end{document}